\documentclass[12pt]{amsart}

\pdfoutput=1

\usepackage{amssymb,latexsym}
\usepackage{enumerate}

\makeatletter
\@namedef{subjclassname@2010}{%
  \textup{2010} Mathematics Subject Classification}
\makeatother

\newtheorem{thm}{Theorem}[section]
\newtheorem{prop}[thm]{Proposition}
\newtheorem{cor}[thm]{Corollary}
\newtheorem{lem}[thm]{Lemma}

\theoremstyle{definition}

\newtheorem{exa}[thm]{Example}

\numberwithin{equation}{section}

\newcommand{\R}{\mathbb R}
\newcommand{\Z}{\mathbb Z}
\newcommand{\N}{\mathbb N}

\newcommand{\al}{\alpha}

\renewcommand{\epsilon}{\varepsilon}

\begin{document}


\baselineskip=17pt



\title[Diophantine approximation]{Inhomogeneous Diophantine
approximation with general error functions }

\author[L. Liao]{Lingmin Liao}
\address{Lingmin Liao\\LAMA UMR 8050, CNRS,
Universit\'e Paris-Est Cr\'eteil, 61 Avenue du
G\'en\'eral de Gaulle, 94010 Cr\'eteil Cedex, France}
\email{lingmin.liao@u-pec.fr}

\author[M. Rams]{Micha\l\ Rams}
\address{Micha\l\ Rams\\Institute of Mathematics\\ Polish
Academy of Sciences\\ ul.
\'Sniadeckich 8, 00-956 Warszawa\\ Poland }
\email{M.Rams@impan.gov.pl}

\date{}


\begin{abstract}
Let $\al$ be an irrational and $\varphi: \N \rightarrow \R^+$ be
a function decreasing to zero.
For any $\al$ with a given Diophantine type, we show some sharp
estimations for the Hausdorff dimension of the set
\[
E_{\varphi}(\al):=\{y\in \R: \|n\al -y\| < \varphi(n) \text{ for
infinitely many } n\},
\]
where $\|\cdot\|$ denotes the distance to the nearest integer.
\end{abstract}

\subjclass[2010]{Primary 28A80; Secondary 37E05, 28A78}

\keywords{Inhomogeneous Diophantine approximation, Hausdorff dimension}

\maketitle

\section{Introduction}

Let $\al$ be an irrational real number. Denote by $\|\cdot\|$
the distance to the nearest integer. A famous result of
Minkowski (\cite{Mink}) in 1907 showed that if $y \not\in \Z+\al
\Z$, then
for infinitely many $n\in\Z$, we have
\[  \| n\al -y \| < \frac{1}{4|n|}.\]
If $n$ is restricted to positive integers only, Khintchine (\cite{Khintchine26}) in 1926 proved that  
for any real number $y$, there exist infinitely many $n\in\N$ satisfying the Diophantine inequalities:
\begin{equation}\label{Dio-ineq}  \| n\al -y \| < \frac{1}{\sqrt{5}n}.\end{equation}
We shall always restrict $n$ to positive integers. Khintchine's resault is equivalent to say that the set 
\[
E(\al, c):=\left\{y\in \R: \|n\al -y\| <
\frac{c}{n} \ \text{ for infinitely many } n
\right\},
\]
is the whole space $\R$ when the constant $c$ equals to $1/\sqrt{5}$.
It is showed by Cassels \cite{Cassels50} in 1950 that the set $E(\al, c)$ is of full measure for any constant $c>0$.

However, if the error function (the right-hand side of (\ref{Dio-ineq})) of the above Diophantine inequalities is replaced by a function decreasing to zero faster than $c/n$, the sizes of the sets in question would be of zero Lebesgue measure and then the Hausdorff dimension is involved.

Define the Diophantine type $\beta(\alpha)$  of $\alpha$ by
\[
\beta(\al):= \sup \{\theta \geq 1: \liminf_{n\to
\infty}n^{\theta} \|n\al\|=0\}.
\]
In 1999, Bernik and Dodson \cite{BD99} proved that the Hausdorff
dimension, denoted by $\dim_H$, of the set \[
E_{\gamma}(\al):=\left\{y\in \R: \|n\al -y\| <
\frac{1}{n^{\gamma}} \ \text{ for infinitely many } n
\right\}\quad (\gamma\geq 1),
\]
 satisfies
\[ \frac{1}{\beta(\alpha)\cdot \gamma}\leq \dim_HE_{\gamma}(\al)
\leq \frac{1}{\gamma}.\]
In 2003, Bugeaud \cite{Bugeaud03}, and independently Schmeling
and Troubetzkoy \cite{TS03} improved the above result. They
showed that for any irrational $\al$,
\[  \dim_HE_{\gamma}(\al)  = \frac{1}{\gamma}.\]

Now let $\varphi: \N \rightarrow \R^+$ be a function decreasing
to zero. 
Consider the set
\[
E_{\varphi}(\al):=\{y\in \R: \|n\al -y\| < \varphi(n) \text{ for
infinitely many } n\}.
\]
This is the set of well-approximated numbers with a general
error function $\varphi$.
It easily follows from the Borel-Cantelli lemma that the Lebesgue measure of $E_{\varphi}(\al)$ is zero whenever the series 
$\sum_{n=1}^{\infty} \varphi(n)$ converges. But on the other hand, it seems hard to obtain a lower bound of the Lebesgue measure of $E_{\varphi}(\al)\cap[0,1]$ when
the series $\sum_{n=1}^{\infty} \varphi(n)$ diverges. For the results on the Lebesgue measure, we refer the readers to \cite{Kurzweil55}, \cite{LN12}, \cite{Kim12}, and the references therein. 

In this paper, we are concerned with the Hausdorff dimension of the set $E_{\varphi}(\al)$. We can find a natural upper bound:
\[ \dim_HE_{\varphi}(\al) \leq \limsup_{n\to\infty}\frac{\log
n}{-\log \varphi(n)}.\]
It can also be proved that for almost all real numbers $\al$,
the above inequality becomes an equality.
However, in \cite{FW06}, Fan and Wu constructed an example which
shows that the equality is not always true. In fact, they found
a Liouville number $\alpha$ and constructed an error function
$\varphi$ such that
\[ \dim_HE_{\varphi}(\al) = \liminf_{n\to\infty}\frac{\log
n}{-\log \varphi(n)}<\limsup_{n\to\infty}\frac{\log n}{-\log
\varphi(n)}.\]
So in general case, the dimension formula seems mystery.

Recently, Xu \cite{Xu} made a progress, he proved the following
theorem.
\begin{thm}[Xu]\label{Xu}

For any $\al$, we have the following estimation
\[
\limsup_{n\to\infty}\frac{\log q_n}{-\log \varphi(q_n)} \leq
\dim_H(E_{\varphi}(\al))\leq \limsup_{n\to\infty}\frac{\log
n}{-\log \varphi(n)},
\]
where $q_n$ denotes the denominator of the $n$-th convergent of
the continued fraction of $\al$.
\end{thm}
As a corollary, Xu proved that for any irrational number
$\al$ with Diophantine type $1$,
\[
\dim_H(E_{\varphi}(\al))= \limsup_{n\to\infty}\frac{\log
n}{-\log \varphi(n)}.
\]

For the simplicity, let us denote \[u_{\varphi}:=
\limsup_{n\to\infty}\frac{\log n}{-\log \varphi(n)} \qquad
l_{\varphi}:=\liminf_{n\to\infty}\frac{\log n}{-\log
\varphi(n)}. \]
In this paper, we prove the following results.
\begin{thm}\label{th-main}
  For any $\al$ with Diophantine type $\beta$, we have
\begin{eqnarray*}
\min \left\{u_{\varphi}, \ \max\left\{l_{\varphi},
\frac{1+u_{\varphi}}{1+\beta}\right\}\right\} \leq
\dim_H(E_{{\varphi}}(\al))\leq u_{\varphi}.
\end{eqnarray*}
\end{thm}

\begin{cor}\label{Cor}

If $\beta\leq 1/u_{\varphi}$, then
\begin{eqnarray*}
\dim_H(E_{{\varphi}}(\al)) =u_{\varphi}.
\end{eqnarray*}
\end{cor}

\begin{exa}\label{ex}
Take $\beta=2$, $u=1/2$ and $l=1/3$. We can construct an
irrational $\alpha$ such that for all $n$, $q_n^2\leq q_{n+1}
\leq 2 q_n^2$. Define \[\varphi(n)=\max\big\{n^{-1/l},
q_k^{-1/l}\big\} \ \ \text{ if } \ q_{k-1}^{u/l}< n \leq
q_k^{u/l}.\]
Then by Corollary \ref{Cor}, we have \[ \lim_{n\to
\infty}\frac{\log q_n}{-\log
\varphi(q_n)}=l<u=\dim_H(E_{\varphi}(\al)).\]
Thus the lower bound of Xu (Theorem \ref{Xu}) is not optimal.

\end{exa}

The next two theorems show that the estimations in Theorem
\ref{th-main} are sharp.
\begin{thm}\label{sharp}
For any irrational $\al$ and for any $0\leq l<u\leq1$, with
$u>1/\beta$, there exists a decreasing function $\varphi : \N
\rightarrow \R^+$, with $l_\varphi=l$ and
$u_\varphi=u$, such that
\[
\dim_H(E_{\varphi}(\al)) =\max\left\{l,
\frac{1+u}{1+\beta}\right\}<u.
\]
\end{thm}

\begin{thm}\label{sharp-upper}
Suppose $0\leq l<u\leq1$. There exists a decreasing function
$\varphi : \N \rightarrow \R^+$, with $l_\varphi=l$ and
$u_\varphi=u$, such that for any $\al$ with $\beta< \infty$,
\[
\dim_H(E_{\varphi}(\al)) =u.
\]
\end{thm}

\bigskip
\section{Three steps dimension}

The goal of this section is to prove Proposition \ref{tech}
which will be the base of our dimension estimation (compare
\cite[Section 3]{Xu}).

As a direct corollary of Proposition
\ref{tech}, we will also give a new proof of Xu's theorem
(Theorem \ref{Xu}) at the end of this section.

Let us start with a technical lemma.
\begin{lem} \label{lem:conv}
Let $1>a>b$ and $1>c>d$. Then for any $\delta \in [0,1]$ we have\[
\frac {\log (\delta a + (1-\delta) c)} {\log (\delta b +
(1-\delta) d)} \geq \min \left(\frac {\log a} {\log b}, \frac
{\log c} {\log d} \right).
\]
\begin{proof}
Denote
\[
s:=\min \left(\frac {\log a} {\log b}, \frac
{\log c} {\log d} \right).
\]
%
Then
\[
\frac {\log (\delta a + (1-\delta) c)} {\log (\delta b +
(1-\delta) d)} \geq\frac {\log (\delta b^s + (1-\delta) d^s)} {\log
(\delta b + (1-\delta) d)}.
\]
By concavity of the function $x\to x^s$, we have
\[
\delta b^s + (1-\delta) d^s \leq (\delta b + (1-\delta) d)^s
\]
and the assertion follows.
\end{proof}
\end{lem}

Let $\alpha$ be an irrational number with Diophantine type $\beta(\alpha)>1$. Recall that $q_n$ is the
denominator of the $n$-th convergent of the continued fraction
of $\al$.
Let $B\geq 1$ and suppose there exists a sequence of natural
numbers $\{n_i\}$ such that
\begin{equation}\label{lim-B}
\frac {\log q_{n_i+1}} {\log {q_{n_i}}} \to B.
\end{equation}
Let $\{m_i\}$ be a sequence of natural numbers
such that $q_{n_i}<m_i\leq q_{n_i+1}$.
  By passing to subsequences, we suppose the limit
\[
N:= \lim_{i\to \infty}\frac {\log m_i} {\log q_{n_i}}
\]
exists. 
Then obviously, $1 \leq N \leq B$.

Let $K>1$.
Denote
\[
E_i := \Big\{y\in \R: ||n\alpha - y|| < \frac 12 q_{n_i}^{-K} \
\text{ for some}\ n\in (m_{i-1}, m_i]\Big\}.
\]
Let
\[
E := \bigcap_{i=1}^\infty E_i \qquad
\text{and}
\qquad
F := \bigcap_{j=1}^\infty \bigcup_{i=j}^\infty E_i.
\]

\begin{prop} \label{tech}
If $\{n_i\}$ is increasing sufficiently fast then
\[
\dim_H E = \dim_H F = S,\]
where
\[
S=S(N,B,K):=\min\left(\frac N K, \ \max \left(\frac 1K, \frac 1
{1+B-N}\right)\right).
\]

\begin{proof}
As $F\supset E$, we only need to get the lower bound for $\dim_H
E$ and the upper bound for $\dim_H F$. For the former, we will use the Frostman Lemma, and for the latter, we will use a natural
cover.

We will distinguish two cases: $B\geq K$ and $B<K$. Notice the following fact. 

{\bf Fact}:
If $B\geq K$ then
\[
\frac N K > \frac 1 {1+B-N}, \quad \text{and} \quad S=\max \left(\frac 1K, \frac 1
{1+B-N}\right).
\]If $B<K$, then
\[
\frac 1 K < \frac 1 {1+B-N}, \quad \text{and} \quad S=\min\left(\frac N K, \ \frac 1
{1+B-N}\right).
\]

Indeed, the second statement follows by noting $1/K < 1/B$. For the first statement, if $N\geq K$ then it is obviously true because the right
hand side is smaller than $1$. Otherwise, we have
\[\frac{K-N}{N}<{K-N},\]
hence
\[
\frac K N <1+  K-N.
\]
Since $B\geq K$, we have
\[
1+B-N  \geq 1+K-N>K/N.
\]


\noindent{\it Distribution of the points}.\\
Now, let us study the distribution of the points $\{n\alpha \
(\text{mod} \ 1)\}$.
Let $\{n_i\}$ be a fast increasing sequence satisfying
(\ref{lim-B}). By passing to a subsequence, we can always assume
that $\{n_i\}$ grows as fast as we wish; the exact conditions on
the rate of growth will be clear later. Denote
\[N_i := m_i - m_{i-1}.\]
By passing to a subsequence, we can suppose that $N_i\geq q_{n_i}$.

The three steps theorem tells us how the points $\{n\alpha \
(\text{mod} \ 1)\}_{n=m_{i-1} +1}^{m_i}$ are distributed on the unit
circle: there are $q_{n_i}$ groups of points, each consisting of
$\lfloor N_i/q_{n_i}\rfloor$ ($\lfloor \cdot \rfloor$ denotes the integer part) points, the distances between points inside each
group are equal to $\xi_i:=\|q_{n_i}\alpha\|$ and the distances
between groups are $\zeta_i:=\|q_{n_i-1}\alpha\|-(\lfloor N_i/q_{n_i}\rfloor-1)\|q_{n_i}\alpha\|$.

In the first case, i.e., $B \geq K$, we have $\xi_i \leq q_{n_i}^{-K}$ for all $i$ big enough,
hence the intervals $[n\alpha - q_{n_i}^{-K}/2,
n\alpha+q_{n_i}^{-K}/2]$ intersect each other (inside each
group). So $E_i$ consists of $M_i:=q_{n_i}$ intervals of length
$y_i:=(\lfloor N_i /q_{n_i}\rfloor-1)\xi_i + q_{n_i}^{-K}$. 
By noting that $\|q_n\alpha\|$ is comparable with $q_{n+1}^{-1}$, we have
\[
y_i = (\lfloor N_i/q_{n_i}\rfloor-1) \xi_i + q_{n_i}^{-K} = q_{n_i}^{-\min(K, 1+B-N)+O(\varepsilon)}.
\]

In the second case, i.e., $B < K$, for big $i$, $E_i$ consists of $N_i$ intervals of
length $z_i:=q_{n_i}^{-K}$.

As $q_{n_{i+1}}\gg q_{n_i+1}$, we can freely assume that for any $\epsilon>0$, each
component of $E_i$ contains at least $M_{i+1}^{1-\varepsilon}$
(in the first case) or $N_{i+1}^{1-\varepsilon}$ (in the second
case) components of $E_{i+1}$. 

\noindent{\it Calculations}.\\
We will distribute a probability
measure $\mu$ in the most natural way: the measure attributed to
each component of $F_i=E_1\cap\ldots\cap E_i$ is the same. Here
we distribute the measure only on those components of $F_i$ that
are components of $E_i$, i.e., at all stages we count only
components completely contained in previous generation sets.

%
%

%

{\bf Case 1}:  $B \geq K$.   At level $i$ we have at least
$M_i^{1-\varepsilon}$ components of $F_i$, each of size
$y_i$ and inside each component of $F_{i-1}$, the components of
$F_i$ are in equal distance $c_i:=\zeta_i-q_{n_i}^{-K}$.

Let $x\in E$. For $y_i \leq r < y_{i-1}$, consider
\begin{equation} \label{eq:fr}
f(r) = \frac {\log \mu(B_r(x))} {\log r}.
\end{equation}
Notice that the convex hull of components of $F_i$ intersecting
$B_r(x)$ has measure at most $3\mu(B_r(x))$ and length at most
$6r$. For simplicity, we can assume that the interval $B_r(x)$
is a convex hull of some components of $F_i$ contained in one
component of $F_{i-1}$. Hence,
\begin{equation} \label{eq:lin1}
f(ny_i + (n-1)c_i) \geq \frac {\log (n M_i^{-(1-\varepsilon)})}
{\log (ny_i + (n-1)c_i)}.
\end{equation}

As the right hand side of equation \eqref{eq:lin1} is the ratio
of logarithms of two functions, both linear in $n$ and smaller
than 1, by Lemma \ref{lem:conv} the minimum of $f(r)$ in range
$(y_i, y_{i-1})$ is achieved at one of endpoints. We have
\begin{equation} \label{eq:fyi}
f(y_i) \geq (1-\varepsilon) \frac {-\log M_i} {\log y_i} = \max \left(\frac 1 K ,\frac
1 {1+B-N}\right) + O(\varepsilon)
\end{equation}
and the same holds for $f(y_{i-1})$. Recalling the fact at the beginning of the proof, we get the lower bound by Frostman Lemma.

The upper bound is simpler: for any $i$, $F$ is contained in
$\bigcup_{n>i} E_n$. Hence, we can use the components of all
$E_n, n>i$ as a cover for $F$. For any $s$ the sum of $s$-th
powers of diameters of components of $E_n$ is bounded by $M_n
y_n^s$, and for $s> \max(\frac 1K,\frac 1 {1+B-N}) + O(\varepsilon)$ it is
exponentially decreasing with $n$. The upper bound then follows by the definition of Hausdorff dimension.

{\bf Case 2}:  $B < K$.  Once again to obtain the
lower bound we will consider the function $f(r)$ given by
\eqref{eq:fr}. However, in this case the components of $F_i$ are
not uniformly distributed inside a component of $F_{i-1}$ but
they are in groups. There are at least $s_i$ groups in distance
$c_i$ from each other, each group is of size $y_i$ and contains
at least $N_i^{1-\varepsilon}$ components. Inside each group
the components of size $z_i$ are in distance $d_i :=
\xi_i-q_{n_i}^{-K}$ from each other.

We need to consider $z_i \leq r < z_{i-1}$. This range can be
divided into two subranges. The equation \eqref{eq:lin1} works
for $y_i \leq r < z_{i-1}$, while for $z_i \leq r < y_i$
the same reasoning gives
\begin{equation} \label{eq:lin2}
f(nz_i + (n-1)d_i) \geq \frac {\log (n N_i^{-(1-\varepsilon)})}
{\log (nz_i + (n-1)d_i)}.
\end{equation}

Like in the first case, Lemma \ref{lem:conv} implies that the
minimum of $f(r)$ in each subrange is achieved at one of
endpoints. We have
\[
f(z_i) \geq (1-\varepsilon) \frac {-\log N_i} {\log z_i} = \frac
N K + O(\varepsilon)
\]
and the same for $f(z_{i-1})$, while $f(y_i)$ is still given by
\eqref{eq:fyi}. Together with the fact at the beginning of the proof,  this gives the lower bound.

To get the upper bound for the dimension of $F$ we can use two
covers. One is given by using the convex hulls of groups of
components of $F_n$ with $n>i$. As in the first case (taking into account the fact that $1/K<1/(1+B-N)$), this cover gives
\[
\dim_H F \leq \frac 1 {1+B-N} + O(\varepsilon).
\]
The other cover consists of
components of $E_n$ with $n>i$. For any $s$ the sum of $s$-th powers of
diameters of components of $E_n$ is bounded by $N_n z_n^s$, and
for $s> \frac N K + O(\varepsilon)$ it is exponentially
decreasing with $n$. We will choose one of the two covers
that gives us the smaller Hausdorff dimension.
\end{proof}
\end{prop}


The statement of Proposition \ref{tech} could be also written in the following way, fixing $B$ and $N$ and varying $K$:
\[
S(N,B,K)= \left\{\begin{array}{ll}
1/K & K < 1+B-N\\
1/(1+B-N) & 1+B-N \leq K \leq N(1+B-N) \\
N/K & K > N(1+B-N).
\end{array}
\right.
\]

\smallskip
By Proposition \ref{tech}, we can directly deduce Theorem
\ref{Xu}.

{\it A new proof of Theorem \ref{Xu}}:\\
The upper bound is easy, and we only show the lower bound. We will
apply Proposition \ref{tech}.
Let $q_{n_i}$ be a sparse subsequence such that
\[
\lim_{i\to \infty} \frac{\log q_{n_i+1}}{-\log
\varphi(q_{n_i+1})}=\limsup_{n\to\infty} \frac{\log q_n}{-\log
\varphi(q_n)}=:L.
\]
By passing to a subsequence, suppose the limit
\[
\lim_{i\to \infty} \frac{\log q_{n_i+1}}{\log q_{n_i}}=:B
\]
exists.
Take $m_i=q_{n_i+1}$. Then $ 1\leq N=B$. Take
$K=NL^{-1}$ and construct the sets $E_i$
and $E$ as in Proposition \ref{tech}.
We can easily check that $E$ is a subset of
$E_{\varphi}(\alpha)$. By Proposition \ref{tech}, we have
\[
\dim_H E_{\varphi}(\alpha) \geq \min\big\{L, \ 1\big\}.
\]
Then the result follows.

%

\bigskip
\section{Proof of Theorem \ref{th-main}}

The upper bound of Theorem \ref{th-main} is trivial by using the
natural covering, hence we will only concern ourselves with the
lower bound.

Note that the lower bound in Theorem \ref{th-main} can be
written as
\[\max \left\{l_{\varphi}, \ \min\left\{u_{\varphi},
\frac{1+u_{\varphi}}{1+\beta}\right\}\right\}.\]
By the result of Bugeaud \cite{Bugeaud03} and Schmeling and Troubetzkoy \cite{TS03}, the Hausdorff dimension of $E_\varphi$ is at least $l_\varphi$.
So we just need to show it
is not smaller than $\min(u_\varphi, (1+u_\varphi)/(1+\beta))$.

We shall suppose that $l_{\varphi}>0$, the case $l_{\varphi}=0$
can be done by a limit argument.
Since the result is known if $l_{\varphi}=u_{\varphi}$, we also
suppose that $l_{\varphi}<u_{\varphi}$.

The Diophantine type of the irrational number $\al$ can be
defined alternatively by
\[\beta=\limsup_{n\to\infty}\frac{\log q_{n+1}}{\log q_n}.\]
Choose a sequence $m_i$ of natural numbers such that
\[
\lim_{i\to\infty} \frac {\log {m}_i} {-\log \varphi({m}_i)} =
u_\varphi.
\]
Let $n_i$ be such that $q_{n_i} < m_i \leq q_{n_i+1}$. By
passing to a subsequence we can assume that
\begin{itemize}
\item[--] the sequence $\log q_{n_i+1}/\log q_{n_i}$ has some
limit $B\in [1, \beta]$,
\item[--] the sequence $\log m_i/\log q_{n_i}$ has some limit
$N\in [1, B]$,
\item[--] the sequence $\{n_i\}$ grows fast enough for
Proposition \ref{tech}.
\end{itemize}
Moreover, we can freely assume that $N>1$: otherwise, by the monotonicity of $\varphi$,  we would have
\[
\lim_{i\to\infty} \frac {\log q_{n_i}} {-\log \varphi(q_{n_i})}
= u_\varphi
\]
and the assertion would follow from Theorem \ref{Xu}.

Take
$K=N/u_\varphi$.
By the definition of $m_i$, for any small $\delta>0$, we have
for all large $i$
\[
\varphi(m_i)\geq (m_i)^{- 1/u_\varphi-\delta} \geq
q_{n_i}^{-K}.
\]
Thus by monotonicity of $\varphi$,
\begin{equation}\label{eq-inclu}
\varphi(n)  \geq q_{n_i}^{-K} \quad \forall n\leq m_i.
\end{equation}



The assumptions of Proposition \ref{tech} are
satisfied, so we can calculate the Hausdorff dimension of the
set $E$ defined in the previous section. By (\ref{eq-inclu}),
$E\subset E_\varphi$, so this gives the lower bound for the
Hausdorff dimension of $E_\varphi$:
\[
\dim_H E_\varphi \geq M(N,B):= \min \left(u_\varphi, \max\left(\frac {u_\varphi} N, \frac 1 {1+B-N} \right)
\right)
\]
and we want to estimate the minimal value of $M$ for $B\in
[1,\beta], N\in [1,B]$.

First thing to note is that increasing $B$ not only decreases
$M(B,N)$ for a fixed $N$ but also increases the range of
possible $N$'s. Hence, the minimum of $M(N,B)$ is achieved for
$B=\beta$. Denote $M(N)=M(N, \beta)$.

We are then left with a simple optimization problem of a
function of one variable. We can write
\[
M(N) = \min\left(u_\varphi, \max\left(\frac {u_\varphi} N, \frac
1 {1+\beta-N}\right)\right).
\]
If $\beta u_\varphi \leq 1$ then $u_\varphi \leq 1/(1+\beta-N)$
for all $N$, hence
\[
\min_N M(N)=u_\varphi \leq \frac {1+u_\varphi} {1+\beta}. \]
Otherwise, as $u_\varphi/N$ is a decreasing and $1/(1+\beta-N)$
an increasing function of $N$, the global minimum over $N$ of
the maximum of the two is achieved at the point $N_0$ where they
are equal: $u_\varphi/N_0 = 1/(1+\beta-N_0)$, that is for
\[
N_0=\frac {u_\varphi(1+\beta)} {1+u_\varphi}.
\]
As $\beta u_\varphi > 1$ implies $1<N_0< \beta u_\varphi \leq
\beta$, $N_0$ is inside the interval $[1,\beta]$, hence this
global minimum is the local minimum we are looking for. Thus,
in this case
\[
\min_N M(N) = M(N_0) = \frac {1+u_\varphi} {1+\beta} <
u_\varphi.
\]
We are done.

\bigskip
\section{Proof of Theorems \ref{sharp} and \ref{sharp-upper}}

{\bf Proof of Theorem \ref{sharp}}: Let $\al$ be of Diophantine
type $\beta>1/u$. Let $q_{n_i}$ be a sparse subsequence of
denominators of convergents such that
\[\beta=\lim_{i\to\infty}\frac{\log q_{n_i+1}}{\log q_{n_i}}.\]
For any $0\leq l<u\leq 1$, define
\[
z=\max\left(l, \frac {1+u} {1+\beta}\right).
\]
Note that $z\leq u$.

Define also a function $\varphi : \N \rightarrow \R$ as follows:\[
\varphi(n):= \max\{n^{-1/l}, \ k_{n_i}^{-1/u} \}, \quad
\text{if} \ k_{n_{i-1}}< n\leq k_{n_i},
\]
where
\[
k_{n_i} = q_{n_i}^{u/z}.
\]

Let $D_1$ be the set
\[
 \{y\in \R: {\text{ for infinitely many}\ i},\ ||n\alpha-y||<k_{n_i}^{-\frac{1}{u}} {\text{ for some}}\ n\in (k_{n_{i-1}}, k_{n_i}]\}
\]
and $D_2$ be the set
\[
\{y\in \R: ||n\alpha-y||<n^{-\frac{1}{l}} \text{ for infinitely
many}\ n\}.
\]

Clearly, $E_\varphi(\alpha) = D_1\cup D_2$.
The Hausdorff dimension of $D_1$ is given by Proposition \ref{tech} (with $B=\beta, K=1/z, N=u/z$):
\[
\dim_H D_1 = \min\left(u, \ \max \left( z, \frac z {(1+\beta)z - u}\right)\right)=z
\]
(the equality is valid both when $z=l$ and $z=(1+u)/(1+\beta)$).

By \cite{Bugeaud03} and \cite{TS03} we have
\[
\dim_H (D_2) = l.
\]
Then the proof is completed.

\medskip
{\bf Proof of Theorem \ref{sharp-upper}}: Construct a sequence
$\{n_i\}_{i\geq 1}$ by recurrence:
\[n_1=2, \quad n_{i+1}=2^{n_i} \ (i\geq 1).\]
Define a function $\varphi : \N \rightarrow \R$ as
$\varphi(n)=n_i^{-1/l}$ for $n\in (n_i, n_i^{u/l})$ and
$\varphi(n)=n^{-1/u}$ elsewhere.

Suppose that $\dim_H(E_{\varphi}(\al))<u$. By Theorem \ref{Xu},
no $q_{m}$ could be between $n_i$ and $n_{i+1}^{l/u}$. Since
$n_i$ go to infinity very fast, $\alpha$ cannot be of finite
type.

\subsection*{Acknowledgements}
M.R. was partially supported by the MNiSW (grant N201 607640, Poland). L.L. was partially supported by the ANR (grant 12R03191A -MUTADIS, France). A part of this paper was written during the visit of L.L. to the NCTS in Taiwan.

\end{document}